\documentclass[letterpaper,10pt]{article}

\usepackage{amssymb}
\usepackage{amsmath}
\usepackage{latexsym}
\usepackage[english]{babel}
\usepackage[draft]{hyperref}

\newcommand{\beq}{\begin{equation}}
\newcommand{\eeq}{\end{equation}}
\newcommand{\beqn}{\begin{eqnarray}}
\newcommand{\eeqn}{\end{eqnarray}}
\newcommand{\p}{\psi}
\newcommand{\var}{\varepsilon}

\newcommand\m{\mu}
\newcommand\et{\widetilde{\varepsilon}}

\newcommand\bbR{\mathbb{R}}
\newcommand\bbC{\mathbb{C}}
\newcommand\bbN{\mathbb{N}}
\newcommand\wk{\widetilde{\kappa}}
\newcommand\ws{\widetilde{\sigma}}
\newcommand\wt{\widetilde{\theta}}
\newcommand\tit{\tilde{\theta}}
\newcommand{\ds}{\displaystyle}

\def\QEDclosed{\mbox{\rule[0pt]{1.3ex}{1.3ex}}} 
\def\qed{{{\hspace{.9\textwidth}\QEDclosed}}\newline} 

\newtheorem{theorem}{Theorem}
\newtheorem{lemma}{Lemma}
\newtheorem{remark}{Remark}

\title{Constructive solution of a bilinear optimal control problem for a Schr\"odinger equation}

\date{}

\author{Lucie Baudouin\footnote{e-mail: {\tt baudouin@laas.fr}}\\
{\it\footnotesize LAAS - CNRS; Universit\'e de Toulouse; 7, avenue du Colonel Roche, F-31077 Toulouse, France.}\\ 
Julien Salomon\footnote{e-mail: {\tt julien.salomon@dauphine.fr}}\\
{\it\footnotesize  CEREMADE,\ Universit\'e Paris-Dauphine, Pl. du M$^{al}$ Lattre de Tassigny, F-75775 Paris, France.}}

\begin{document}

\maketitle

\begin{abstract}
Often considered in numerical simulations related to the control of quantum systems, the so-called monotonic schemes
have not been so far much studied from the functional analysis point
of view. Yet, these procedures provide an efficient constructive
method for 
solving a certain class of optimal control problems. This paper aims
both at extending the results already available about these algorithms
in the finite dimensional case ({\it i.e.}, the time-discretized case)
and at completing those of the continuous case. This paper starts with
some results 
about the regularity of a functional related to a wide class of model in quantum chemistry. 
Those enable us to extend an inequality due to \L ojasiewicz to the infinite
dimensional case. Finally, some inequalities proving the Cauchy character of
the monotonic sequence are obtained, followed by an estimation of the rate of convergence. 
\end{abstract}



\section{Introduction}

Following the increasing interest of the chemists community for optimal 
control of quantum systems 
\cite{rabitz2,shi} and the successful laboratory
demonstration of control over molecular phenomena (see, {\it
  e.g.}, \cite{assion,brixner1,weinacht} and more recently \cite{kasparian,vogt}), some mathematical studies of the 
models involved in this topic have  
been carried out, see e.g. \cite{beauchard,lebris}.
In this way, it has been proved in recent papers \cite{baudouin2,pilot} that a
wide class of optimization problems considered by chemists are well
posed. Yet, these proofs are not constructive and
consequently do not give rise to concrete numerical methods to approximate their
solutions.

On the other hand, at numerical simulation level \cite{brown,rabitz3}, various kind of procedures exist and
show a good efficiency. Among them, the so-called monotonic algorithms have
demonstrated their efficiency on several problems. In a recent paper, a
study of the time-discretized algorithms \cite{salomon-cv} have been presented and first
functional analysis results have been obtained about the continuous case \cite{I-K,salomon-cdc05}. 

The aim of this paper is to complete these works by providing
 general proofs of convergence of the optimizing
 sequences. Consequently, we obtain a constructive method, independent
 of time or space discretization to compute  
critical points (and sometimes extrema, see Remark \ref{rem:conv_mini}) of the cost functional under consideration.\\
\indent Let us briefly present the monotonic schemes in the simple
case of ordinary differential equations (ODE). Let $A,B,C$ be three square
matrices in $\mathcal M_n(\bbR)$, $C$ being symmetric positive, $\alpha>0$ and $T>0$. 
Consider  the optimal control problem corresponding to the maximization of the functional $J$ defined by:
\beq\nonumber
J(v)=y(T)\cdot C y(T) -\alpha\int_0^T v^2(t)dt,
\eeq
where $''\cdot''$ denotes the usual scalar product of $\bbR^n$. Here, the state
$y:[0,T]\rightarrow \bbR^n$ and the control $v:[0,T]\rightarrow \bbR$ are linked
by the ODE:
\beq\nonumber
\left\{
\begin{array}{l}
y'(t)=\big(A+v(t)B\big)y(t),~\forall t\in(0,T)\\
y(0)=y_0
\end{array}
\right.
\eeq
the initial condition $y_0$ being fixed. \\
Given two controls $v$ and
$\widetilde{v}$ and the corresponding states $y$ and $\widetilde{y}$, we first note that:
\beqn\nonumber
J(\widetilde{v})-J(v)&=&\big(\widetilde{y}(T)-y(T)\big)\cdot  C
\big(\widetilde{y}(T)-y(T)\big)+2\big(\widetilde{y}(T)-y(T)\big)\cdot  Cy(T)\\\nonumber &&
-~\alpha\int_0^T \big(\widetilde{v}(t)-v(t)\big)(\widetilde{v}(t)+v(t)\big)dt.
\eeqn
We then introduce an auxiliary function $z:[0,T]\rightarrow \bbR^n$ associated
to $y$ and $v$ by
\beq\nonumber
\left\{
\begin{array}{lll}
z'(t)&=&-\big(A^*+v(t)B^*\big)z(t),\\
z(T)&=&Cy(T)
\end{array}
\right.
\eeq
where $A^*$ and $B^*$ are the transposed matrices of $A$ and $B$.\\
Focusing on the second term of the right hand side of this equation, we get:
\beq\nonumber
\big(\widetilde{y}(T)-y(T)\big)\cdot Cy(T)
=\int_0^T\big(\widetilde{v}(t)-v(t)\big) B\widetilde{y}(t)\cdot z(t)dt.
\eeq
Thus, we finally obtain:
\beq\nonumber
J(\widetilde{v})-J(v)\!\!=\!\!\big(\widetilde{y}(T)-y(T)\big)\cdot  C
\big(\widetilde{y}(T)-y(T)\big)
+\alpha\int_0^T\!\!
\big(\widetilde{v}(t)-v(t)\big)\!\left(\frac{2}{\alpha}B\widetilde{y}(t)\cdot z(t) -\widetilde{v}(t)-v(t)\right) \! dt.
\eeq
A simple way to guarantee that $\widetilde{v}$ gives a better cost functional value than $v$, is to impose that:
\beq\label{mono-cond}
\big(\widetilde{v}(t)-v(t)\big)\!\left(\frac{2}{\alpha}B\widetilde{y}(t)\cdot z(t) -\widetilde{v}(t)-v(t)\right)\geq0.
\eeq
Following this approach, the sequence
$(v^k)_{k\in\bbN}$ defined iteratively by the implicit equation 
$v^{k+1}=\frac{1}{\alpha} By^{k+1}(t)\cdot z^{k}(t),$
where $y^{k+1}$ and $z^{k}$ correspond to $v^{k+1}$ and $v^{k}$ respectively,
optimizes $J$ monotonically  since
$$ J(v^{k+1})-J(v^k) = \big(y^{k+1}(T)-y^k(T)\big)\cdot  C
\big(y^{k+1}(T)-y^k(T)\big)
+\alpha\int_0^T \big(v^{k+1}(t)-v^k(t)\big)^2dt\geq 0.$$
In this article, we prove the convergence of generalizations of this algorithm
towards a critical point of $J$ in the case of the Schr\"odinger
partial differential equation:
\beq\nonumber
i\partial_t \p (x,t)-[H-\mu(x)\var(t)]\p (x,t)=0.
\eeq
This equation governs the evolution of a quantum system, described by
its wave function $\p$, that interacts with a laser pulse of amplitude $\var$, the control variable. The
factor $\mu$ is the dipole moment operator of the system. In
what follows, $H=-\Delta+V$ where $\Delta$ is
the Laplacian operator and $V=V(x)$ the electrostatic potential in which
the system evolves. We refer to \cite{rabitz3} for more details about
models involved in quantum control.\\
\indent The paper is organized as follows: we start in Section~\ref{sec2} with some
necessary results about the linear and nonlinear Schr\"odinger equations
involved in the problem we are considering. 
We then present the optimization problem in  
Section~\ref{sec3}, and claim some
regularity results about the corresponding cost functional in
Section~\ref{sec4}. We introduce in Section~\ref{sec5} an important tool for
proving the convergence of the sequence, 
namely the \L ojasiewicz inequality and some of its generalizations. The
definition of the monotonically optimizing sequence is given in
Section~\ref{sec6} where some useful properties are also claimed. The
convergence of 
the sequences is proved in Section~\ref{sec7} and a first result about
their rate of convergence follows in the last section.\\
\indent Throughout this paper, $T$ is a positive real number representing the time of
control of a physico-chemical process. We denote by $L^2$ and $L^\infty$  the
spaces $L^2(\bbR^3,\bbC)$ and $L^\infty(\bbR^3,\bbC)$,
$W^{p,\infty}(\bbR^3,\bbR)$ with $p\in 
[1,+\infty)$ by $W^{p,\infty}$, the Sobolev space
$H^2(\bbR^3,\bbC)$ by $H^2$ and $L^{p}(0,T;X)$, with $p\in 
[1,+\infty)$ denotes the usual Lebesgue space taking its values in a Banach
space $X$. 
We also use the notation $\langle~.~|~.~|~.~\rangle$ and $\langle~.~,~.~\rangle$
defined by: 
\beq\nonumber
\langle f |A|g\rangle=\ds\int_{\bbR^3}\overline{f(x)}Ag(x)dx
,\ \langle f,g\rangle=\ds \int_{\bbR^3}\overline{f(x)}g(x)dx,
\eeq
where
$f$ and $g$ are in $L^2$ and $A$ is an operator on $L^2$. To simplify our
notation, the space variable $x$ will often be omitted. Finally, for $h\in
L^p(0,T;X)$, $p\in ]1,\infty]$, we recall that $\|h\|_{L^p(0,T;X)}=\left\|t\mapsto
  \|h(t)\|_{X} \right\|_{L^p(0,T)}$. Finally, we denote by $Im(z)$ and
$Re(z)$ the imaginary and the real part of a complex number $z$.\\

\section{Preliminary existence results}\label{sec2}

The sequences we study in this paper are defined through iterative
resolutions of Schr\"odinger equations. Before
introducing the relevant framework of our study, we present here some
necessary preliminary existence and regularity results concerning these
equations. The first one will correspond later to the initialization step in
the definition of the sequences. This lemma is a corollary of a
general result on time dependent hamiltonians (see \cite{R-S}, p285,
Theorem~X.70) but for the sake of clarity, we present here an approach
using other techniques also useful in the proof of the next lemmas.

\begin{lemma}\label{schrod0}
Let $\mu$ and $V$ belong to  $W^{2,\infty}$
and let  $H=-\Delta+V$.
If  $\varepsilon\in L^{2}(0,T)$ and $\psi_0\in H^2$, the equation
\beq\label{eqschrod0}
\left\{
\begin{array}{l}
i\partial_t \p (x,t)-[H(x)-\mu(x)\var(t)]\p (x,t)=0\\
\p(x,0)=\p_0(x)
\end{array}
\right.
\eeq
has a unique solution $\p\in L^\infty(0,T;H^2)\cap W^{1,\infty}(0,T;L^2)$. Moreover:
\beq\label{CL2N}
\forall t\in [0,T],\quad  \|\p (t)\|_{L^2}=\|\psi_0\|_{L^2}.
\eeq 
\end{lemma}
{\bf Proof:} \indent  One can also read a similar proof in \cite{baudouin} but we give here some details.
It is well known (see \cite{Caze} for instance) that for any $T>0$ and $u_0 \in H^2$, the Schr\"odinger equation
$$
\left\lbrace
\begin{array}{l}
i\partial_t u(x,t)+\Delta u(x,t)=0,\quad x\in\mathbb{R},~t\in [0,T]\\
u(x,0)=u_0(x),\quad x\in\mathbb{R}
\end{array} \right.
$$
has a unique solution $u(t)=S(t)u_0$ such that $u\in C([0,T];H^2)\cap C^1([0,T];L^2)$, where $(S(t))_{t\in {\mathbb R}}$ denotes the free Schr\"odinger semi-group ${\rm e}^{it\Delta }$. Moreover, for all $t\in [0,T]$ we have 

\begin{equation}
\label{Sfree}
\|u(t)\|_{H^2} = \|S(t)u_0\|_{H^2} =   \|u_0\|_{H^2}.
\end{equation}

Let $\lambda >0$ be a given positive number which will
be fixed hereafter and denote by $Y$ the space $C([0,T];H^2 )$ endowed
with the norm $\| \p\|_Y= \sup_{t\in [0,T]}{\rm e}^{-\lambda t}\| \p(t)\|_{H^2}.$
The solution of equation~(\ref{eqschrod0}) is obtained equivalently as a solution to the integral equation
$$ \p(t)=S(t) \p_0 + i \int_0^t S(t-s)W(s) \p(s)\,ds$$
where $W(x,t) = -V(x) + \mu(x)\var(t)$ for all $t\in[0,T]$, $x\in\bbR^3$.
We are going to show that this equation has a unique solution in $Y$, by proving that 
operator $\Phi$ defined by
$$\Phi ( \p)(t) = S(t) \p_0 + i \int_0^t S(t-s)W(s) \p(s)\,ds$$
has a unique fixed point in a closed ball $B_R = \left\{ \p\in
Y \; ; \; \| \p\|_Y \leq R \right\}$ for suitable $R$.\\

If $ \p \in B_R$, then $\| \p(s)\|_{H^2 } \leq {\rm e}^{\lambda s} \| \p\|_Y \leq R{\rm
e}^{\lambda s}$ and since $W\in L^2 (0,T; W^{2,\infty})$, we can set $\rho>0$ such that
$\|W\|_{L^2 (0,T; W^{2,\infty})} \leq \rho$. Using estimate (\ref{Sfree}) and 
Cauchy-Schwarz inequality we obtain
$$
\|\Phi ( \p)(t)\|_{H^2 } \leq  \| \p_0\|_{H^2 } +   \ds\int_0^t\|W(s) \p(s)\|_{H^2 }ds
\leq \| \p_0\|_{H^2 } + \rho~  R \left(\ds\int_0^t {\rm e}^{2\lambda s}ds\right)^{\frac 12}.
$$
It follows that if $R>0$ is large enough so that $ \| \p_0\|_{H^2 } \leq \dfrac R2$ and if we choose $\lambda > 2\rho^2 $, then
$$
\|\Phi ( \p)\|_Y\leq\sup_{t\in [0,T]} {\rm e}^{-\lambda t}\| \p_0\|_{H^2 } + \rho~  R \left(\ds\int_0^t {\rm e}^{2\lambda (s-t)}ds\right)^{\frac 12}
\leq \dfrac R2 + \frac{\rho~ R}{\sqrt{2\lambda}} ~\leq~ R.$$
This means that $\Phi $ maps $B_R$ into itself. Then, for $ \p_1, \p_2 \in
B_R$ it is clear that
\beq\nonumber
\|(\Phi ( \p_1) - \Phi ( \p_2))(t)\|_{H^2 } 
\leq  \ds\int_0^t\|W(s)( \p_1- \p_2)(s)\|_{H^2 } ds 
\leq \rho~  \| \p_1- \p_2\|_Y \left(\frac{{\rm e}^{2\lambda t}-1}{2\lambda }\right)^{\frac 12},
\eeq
and since $\lambda $ has been appropriately chosen, this proves that $\Phi $ is a
strict contraction from $B_R$ into itself as
$$\|(\Phi ( \p_1) - \Phi ( \p_2))\|_Y 
\leq \rho~  \| \p_1- \p_2\|_Y \!\!\sup_{t\in [0,T]}\left(\frac{1-{\rm e}^{-2\lambda t}}{2\lambda }\right)^{\frac 12}\!\!\!\!
\leq \frac{\rho~ }{\sqrt{2\lambda}} \| \p_1- \p_2\|_Y
\leq \frac 12 \| \p_1- \p_2\|_Y $$
and therefore $\Phi $ has a unique fixed point, yielding the solution of equation~(\ref{eqschrod0}) in 
$L^\infty(0,T;H^2)$.
One can notice that uniqueness is not only true in $B_R$ but also easily proved using the norm in  $L^\infty (0,T;L^2)$.
Moreover, calculating $Im
\ds\int_{\mathbb{R}}(\ref{eqschrod0}).\overline{\p}(x)\,dx$, one can prove the
conservation of the $L^2$-norm (\ref{CL2N}) and finally, using equation
(\ref{eqschrod0}), it is easy to obtain that $\p\in W^{1,\infty}(0,T;L^2)$. 

 \qed

We will also have recourse to a similar lemma, dealing with equation (\ref{eqschrod0}) with a non zero source term.

\begin{lemma}\label{schrod1}
Let $H$, $\mu$, $\var$ be as above and $\p \in L^\infty(0,T;H^2)$. Given $\var'\in
L^{2}(0,T)$, the equation:
\beq\label{eqschrod1}
\left\{
\begin{array}{l}
i\partial_t \p' (x,t)-[H(x)-\mu(x)\var(t)]\p' (x,t)=-\mu (x) \var'(t)\p(t,x)\\
\p'(x,0)=0
\end{array}
\right.
\eeq
has a unique solution $\p'\in L^\infty(0,T;H^2)\cap
W^{1,\infty}(0,T;L^2)$. Moreover the following estimate holds:
\beq\label{estimp}
\|\p'\|_{L^\infty(0,T;L^2)}\leq 2\|\mu\|_{L^\infty}\|\var'\|_{L^1(0,T)}\|\p\|_{L^\infty(0,T;L^2)}.
\eeq
\end{lemma}
{\bf Proof:} \indent  The key point to prove the existence of a solution for
(\ref{eqschrod1}) is to underline the fact that the source term  
$f(x,t) = \mu (x) \var'(t)\p(t,x)$ of this linear Schr\"odinger equation
belongs to $ L^2(0,T;H^2)$. It is then very classical to get from the
Lemma~\ref{schrod0} the existence and uniqueness of a solution $\p'$ to
equation~(\ref{eqschrod1}) in $ L^\infty(0,T;H^2)\cap
W^{1,\infty}(0,T;L^2)$.
Consider now $\varphi\in C([0,T])$ defined on $[0,T]$ by  $\varphi (t) =
\|\p'(t)\|^2_{L^2}$. We have:  
\beq
\frac{d}{dt}\varphi(t)
=2Re~\left\langle\p'(t),\frac{H-\mu\var(t)}{i}\p'(t)-\frac{\mu\var'(t)}{i}\p(t)\right\rangle
=-2\var'(t)~Im \left\langle\p'(t)|\mu|\p(t) \right\rangle.\label{cader}
\eeq
Moreover, there exists $t_0$
such that:
 $\varphi(t_0)=\sup_{t\in[0,T]}\left\{ \|\p'(t)\|^2_{L^2} \right\}.$
We suppose $\var'\neq 0$ and $\psi\neq 0$, so that
$\varphi(t_0)\neq 0$ by uniqueness of the
solution of (\ref{eqschrod1}). Since $\p'(x,0)=0$ for all $x\in \bbR^3$, the integration of (\ref{cader}) between $0$ and $t_0$ yields $~\varphi(t_0)=\ds\int_0^{t_0}-2\var'(t)Im\langle\p'(t)|\mu|\p(t)\rangle dt~$, then:
$$
\varphi(t_0)=\|\psi'(t_0)\|^2_{L^2}\leq
\|\mu\|_{L^\infty}\|\p'(t_0)\|_{L^2}\ds\int_0^{T}2|\var'(t)|\|\p(t)\|_{L^2} dt.
$$
Since $\|\psi'(t)\|_{L^2}\leq \|\psi'(t_0)\|_{L^2}$ for all $t\in[0,T]$, we obtain
$$
\| \psi'(t)\|_{L^2}\leq  \|\psi'(t_0)\|_{L^2}
\leq 2 \|\mu\|_{L^\infty}\|\p\|_{L^\infty(0,T;L^2)}\int_0^{T}|\var'(t)| dt,
$$
what ends the proof of estimate (\ref{estimp}).

\qed
Finally, we claim a last result that will be useful to tackle the problems
related to a nonlinear Schr\"odinger equation we encounter in this study. Actually the nonlinearity we consider here is the one that appears naturally in the adjoint system from a quadratic cost functional (as  $J$ is in (\ref{defJ})), even when the state equation is linear.

\begin{lemma}\label{schrod2}
Let $H$, $\mu$, $\var$ and $\psi_0$ defined as above. Given $\chi\in L^\infty(0,T;H^2)$,
the nonlinear Schr\"odinger equation:
\beq\label{eqschrod2}
\left\{
\begin{array}{l}
i\partial_t \p (x,t)-\left[H(x)-\mu(x)\var(t) + Im\left \langle \chi (t)| \mu | \p(t)
\right \rangle \mu(x)\right]\p (x,t)=0\\
\p(x,0)=\p_0(x)
\end{array}
\right.
\eeq
has a unique solution $\p \in L^\infty(0,T;H^2)\cap W^{1,\infty}(0,T;L^2)$. 
\end{lemma}
{\bf Proof:} \indent  - \textit{First Step} -\\
Let $u$ and $\chi\in H^2$, we denote the nonlinear term by 
$~F(u) = Im \langle \chi | \mu | u \rangle \mu u~$
and we can prove that one has the following estimates: $~\exists C=C(\chi,\mu)>0 $ such that  \\
\begin{equation} \label{Ci}
\forall u,v\in L^2,~~\|F(u)-F(v)\|_{L^2} \leq C(\|u\|_{L^2} + \|v\|_{L^2}) \|u-v\|_{L^2}
\end{equation}
\begin{eqnarray}
\label{Cii}\forall u,v\in H^2,~~\|F(u)-F(v)\|_{H^2} &\leq& C(\|u\|_{L^2} + \|v\|_{H^2}) 
\|u-v\|_{H^2}\\
\label{Ciii}
\|F(u)\|_{H^2} &\leq& C \|u\|_{L^2}  \|u\|_{H^2}
\end{eqnarray}

Indeed
\begin{eqnarray*}
\|F(u)-F(v)\|_{L^2} &\leq& \left\| Im \langle \chi | \mu | u \rangle \mu u -  
Im \langle \chi | \mu | v \rangle \mu v  \right\|_{L^2} \\
&\leq& \left\| Im \langle \chi | \mu | u \rangle \mu (u-v) \right\|_{L^2} +  
\left\| Im \langle \chi | \mu | (u-v) \rangle \mu v  \right\|_{L^2}\\
&\leq& \|\mu\|^2_{L^\infty}\|\chi\|_{L^2}(\|u\|_{L^2} + \|v\|_{L^2}) \|u-v\|_{L^2}
\end{eqnarray*}
which proves (\ref{Ci}). Now, we have to establish (\ref{Cii}) and (\ref{Ciii}).
First of all we have
$$ \|F(u)-F(v)\|^2_{H^2} = \|F(u)-F(v)\|^2_{L^2} + \|\Delta F(u)-\Delta F(v)\|^2_{L^2}.$$
The first term of the right hand side is conveniently bounded in (\ref{Ci}). 
Moreover
\begin{eqnarray*}
\|\Delta F(u)-\Delta F(v)\|_{L^2}
& \leq&\left\| Im \langle \chi | \mu | u-v \rangle \Delta(\mu (u-v)) \right\|_{L^2}
+\left\| Im \langle \chi | \mu | v \rangle \Delta(\mu v) \right\|_{L^2} \\
&\leq&  \|\mu\|^2_{W^{2,\infty}}\|\chi\|_{L^2}(\|u\|_{L^2} + \|v\|_{H^2}) \|u-v\|_{H^2}\\
&\leq& C(\|u\|_{L^2} + \|v\|_{H^2}) \|u-v\|_{H^2}.
\end{eqnarray*}
Then, $F$ is locally lipschitz in $H^2$. 
Therefore, taking $v=0$, we also get  (\ref{Ciii}).\\

- \textit {Second Step} -\\
The proof of a local-in-time result is based again on a fixed point
theorem. We begin by fixing an arbitrary time $T>0$ and considering $\tau \in
]0,T]$. We also consider the functional  
$$\xi : \psi\longmapsto U(~.~,0)\psi_0 - i\int_0^.U(~.~,s)F(\psi(s))\,ds,$$
where $\{U(t,s),s,t\in[0,T]\}$ is the propagator associated with the operator
$H -\mu\var$ and induced by Lemma~\ref{schrod0} (such that
$U(t,s) \in \mathcal L(H^2)$ - for details, see \cite{baudouin3}), and the set 
$$ B=\{ v\in L^{\infty}(0,\tau;H^2), \|\psi\|_{L^{\infty}(0,\tau;H^2)}\leq 2M\|\psi_0\|_{H^2} \}.$$
where $M$ satisfies $\forall v\in H^2, \|U(t,s)v\|_{H^2}\leq M \|v\|_{H^2}$.\\

If $\tau>0$ is small enough, the functional $\xi$ maps $B$ into itself  and is a strict contraction in the Banach space $L^{\infty}(0,\tau;H^2)$. 
Indeed, on the one hand, from estimate (\ref{Ciii}), if $\psi\in B$, we have for all  
$ t\in [0,\tau]$:
\begin{eqnarray*}
\|\xi(\psi)(t)\|_{H^2} &\leq& 
\left\|U(t,0)\psi_0 - i\int_0^t U(t,s)F(\psi(s))\,ds \right\|_{H^2}\\
&\leq& M \|\psi_0\|_{H^2}
+ \tau M \|F(\psi)\|_{L^{\infty}(0,\tau;H^2)}\\
&\leq& M \|\psi_0\|_{H^2}
+ \tau C M \|\psi\|_{L^{\infty}(0,\tau;L^2)} \|\psi\|_{L^{\infty}(0,\tau;H^2)}\\
&\leq& M \|\psi_0\|_{H^2} + 4\tau C M^3\|\psi_0\|_{H^2}^2.
\end{eqnarray*}
Then, if we choose $\tau$ such that 
$4\tau C M^2\|\psi_0\|_{H^2} < 1$ we obtain
$\|\xi(\psi)\|_{L^{\infty}(0,\tau;H^2)} \leq  2M \|\psi_0\|_{H^2}$
and $\xi(\psi)$ belongs to  $B$.
On the other hand, if $\psi_1$ and $\psi_2\in B$, then for all $t$ in $[0,\tau]$ we have,
\begin{eqnarray*}
\|\xi(\psi_1)(t)-\xi(\psi_2)(t)\|_{H^2}\!\!\!\! &\!\!=\!\!& \left\|\int_0^t U(t,s)\left(F(\psi_1(s))-F(\psi_2(s))\right)\,ds\right\|_{H^2}\\
&\!\!\leq\!\!& ~C M \left(\|\psi_1\|_{L^{\infty}(0,\tau;L^2)} 
+ \|\psi_2\|_{L^{\infty}(0,\tau;H^2)}\right) \int_0^t  \|\psi_1(s)-\psi_2(s)\|_{H^2}\,ds\\
&\!\!\leq\!\!& ~4\tau C M^2 \|\psi_0\|_{H^2}~\|\psi_1-\psi_2\|_{L^{\infty}(0,\tau;H^2)}
\end{eqnarray*}
with $4\tau C M^2\|\psi_0\|_{H^2} < 1$.
Therefore, from a usual fixed point theorem, we can deduce existence and uniqueness in the set $B$, then in $L^{\infty}(0,\tau;H^2)$, for $\tau > 0$ small enough, of the solution of equation
\begin{equation}
\label{CU}
\psi(t)=U(t,0)\p_0 - i\int_0^tU(t,s)F(\p(s))\,ds 
\end{equation}
which is in fact equivalent to equation~(\ref{eqschrod2}). Moreover, using
(\ref{eqschrod2}), it is easy to prove that $\partial_t \p$ belongs to
$L^{\infty}(0,\tau;L^2)$. 

The last point is then to prove the uniqueness of the solution $u$ of (\ref{eqschrod2}) in the space $L^{\infty}(0,\tau;H^2) \cap W^{1,\infty}(0,\tau;L^2)$.
Let $\p_1$ and $\p_2$ be two solutions of (\ref{eqschrod2}) 
and $w$ equal to $\p_1-\p_2$. 
Then $w(0)=0$ and 
\begin{equation}\label{Cu-v}
i\partial_t w - [H(x)-\mu(x)\var(t)] w  = F(\p_2)-F(\p_1).
\end{equation}
Calculating $Im \ds\int_{\mathbb{R}} (\ref{Cu-v}).\overline{w}(x)\,dx$ and using (\ref{Ci}) we obtain $\dfrac d{dt}(\|w\|^2_{L^2})\leq C \|w\|^2_{L^2}$
and uniqueness follows by Gronwall lemma.
Hence the proof of uniqueness, existence and regularity of the solution of 
equation~(\ref{eqschrod2}) in $\mathbb R\times[0,\tau]$ for any time  
$\tau < \dfrac 1{4C M^2\|\p_0\|_{H^2}}$.\\

- \textit{Third Step} -\\
Now, the goal is to obtain an a priori estimate 
of the solution in $W^{1,\infty}(0,T;L^2)\cap L^\infty(0,T;H^2)$ for any arbitrary time $T$, in order to prove that the local solution we obtained previously exists globally because we have a uniform bound on the norm $\|\p(t)\|_{H^2}+ \|\partial_t \p(t)\|_{L^2}$.

Actually, since equation~(\ref{eqschrod2}) is equivalent to the integral equation (\ref{CU})  and since it is easy to prove the conservation of the $L^2$-norm of the solution, we have, 
\begin{eqnarray*}
\|\p(t)\|_{H^2} &\leq& \|U(t,0)\p_0\|_{H^2}  +\left\|\int_0^tU(t,s)F(\p(s))\,ds\right\|_{H^2} \\
&\leq& M\left\|\p_0\right\|_{H^2} 
+ M C \int_0^t \left\| \p(s)\right\|_{L^2}\left\|\p(s)\right\|_{H^2} \,ds\\
&\leq& C_{0,T} \left( 1+  \int_0^t \left\|\p(s)\right\|_{H^2} \,ds\right)
\end{eqnarray*}
where $C_{0,T}>0$ is a generic constant depending on the time $T$, on $\mu$, $\chi$
and on $\|\p_0\|_{H^2}$. We finally obtain from Gronwall lemma and from 
equation~(\ref{eqschrod2}),  that  $\|\p(t)\|_{H^2}+ \|\partial_t
\p(t)\|_{L^2} \leq C_{0,T}$ for all $t\in [0,T]$. Hence the proof of Lemma~\ref{schrod2}.

\qed
\section{Optimization problem}\label{sec3}
Let us now present the optimization problem we are dealing with in this paper.
Let $O$ be a positive symmetric bounded operator on $H^2$ and
$\alpha$ and $T$ two positive real numbers. Given $\psi_0\in H^2$, we consider the cost
functional $J$ defined on $L^2(0,T)$ by: 
\beq\label{defJ}
J(\varepsilon)=\langle\p(T)|O|\p(T)\rangle-\alpha\int_0^T\varepsilon^2(t)dt,
\eeq 
where $\p$ is the solution of (\ref{eqschrod0}). 
In all the sequel we suppose that $\|\p_0\|_{L^2} = 1$. 
The existence of a minimizer for similar cost functionals (with the opposite sign)
has been obtained in \cite{baudouin}, \cite{baudouin2} and \cite{pilot} and
follows from the construction of a minimizing sequence and a compactness lemma
(Aubin's lemma). Here, the point is to maximize the functional $J$ and as usual, at the maximum of $J$, the Euler-Lagrange critical point equations are satisfied.
A standard way to write these equations is to use
 a Lagrange multiplier  $\chi(x,t)$ usually called
{\it adjoint state}. The following
critical point equations
are thus obtained, for $x\in\bbR^3$ and $t\in (0,T)$:
\beq\label{schrod0eq}
\left\{\begin{array}{l}
i\partial_t \psi (x,t)-[H(x)-\mu(x)\var(t)]\psi (x,t)=0,\\
\psi(x,0)=\p_0 (x),
\end{array}
\right.
\eeq
\beq\label{schrodadj}
\left\{\begin{array}{l}
i\partial_t \chi (x,t)-[H(x)-\mu(x)\var(t)]\chi (x,t)=0,\\
\chi(x,T)=O\p (x,T),
\end{array}
\right.
\eeq

\beq\nonumber
\alpha \var (t) + Im \langle \psi (t) |\mu|\chi(t) \rangle=0.
\eeq
The existence of $\chi\in L^\infty(0,T;H^2)$ results from an adaptation of 
Lemma~\ref{schrod0}, as for $\p(T)\in H^2$ since equation
(\ref{schrod0eq}) is actually equation (\ref{eqschrod0}). In what follows, we
also consider the linearized equation 
of (\ref{schrodadj}):
\beq\label{schrodadj1}
\left\{
\begin{array}{l}
i\partial_t \chi' (x,t)-[H(x)-\mu(x)\var(t)]\chi' (x,t)=-\mu (x) \var'(t)\chi(t,x)\\
\chi'(x,T)=O\p'(T),
\end{array}
\right.
\eeq
where $\var'\in L^2(0,T)$ and $\p'$ is the solution of (\ref{eqschrod1}),
corresponding to the
solution $\p$ of (\ref{schrod0eq}). The existence of $\chi'\in L^\infty(0,T;H^2)$ 
follows from Lemma~\ref{schrod1}. 
The analysis done in the proof of estimate (\ref{estimp}) gives in this
case:
\beq\label{estimadj}
\|\chi'(t)\|_{L^2}\leq 2\|\mu\|_{L^\infty} \|\var'\|_{L^1(0,T)}\|\chi\|_{L^\infty(0,T;L^2)}+\|\chi'(T)\|_{L^2}.
\eeq
Since $\chi(T)=O\p(T)$, $\chi'(T)=O\p'(T)$ and from (\ref{estimp}) and the conservation of the $L^2$-norm, we obtain
\beq\label{estimadj2}
\begin{array}{ccl}
\|\chi'\|_{L^\infty(0,T;L^2)}&\leq&   2\|\mu\|_{L^\infty}
\|\var'\|_{L^1(0,T)}\|\chi\|_{L^\infty(0,T;L^2)}+2 \|O\|_*\|\mu\|_{L^\infty}\|\var'\|_{L^1(0,T)}
\|\p\|_{L^\infty(0,T;L^2)}\\
&\leq&  4  \|\mu\|_{L^\infty}\|O\|_* \|\var'\|_{L^1(0,T)},
\end{array}
\eeq
where $\|O\|_*$ denotes the operator norm
of $O$ on $L^2$. 

\section{Properties of the functional $J$}\label{sec4}
We begin with some properties about the regularity of the cost functional $J$.

\subsection{Gradient of $J$}

We start with some first order properties. As often, the use of the
adjoint state $\chi$ allows us to simplify the computation of the derivative of
$J$. This result is the purpose of the next lemma.   

\begin{lemma}\label{grad}
The cost functional $J$ is differentiable on
$L^2(0,T)$ and its gradient can be expressed by
\beq\label{dJ}
(\nabla J (\var),\var')=-2\int_0^T \big(\alpha \var (t) + Im
\langle\chi(t)|\mu|\psi(t)\rangle\big)\var' (t)dt,
\eeq 
where $(\cdot,\cdot)$ is the usual inner product on $L^2(0,T)$ and $\psi$ and
$\chi$ are the solutions of $(\ref{schrod0eq})$ and $(\ref{schrodadj})$.
\end{lemma}
{\bf Proof:} \indent  We only give here a sketch of the proof. The details can be found
in reference \cite{baudouin} for a slightly different
cost functional. The main point is to prove the differentiability of the
functional  
$\phi :
\var  \in L^2(0,T) \mapsto \p(T),$
 where  $\p$ is the solution of equation (\ref{schrod0eq}).
Actually, one can prove that the solution $\p'$ of (\ref{eqschrod1}) 
is such that 
$D\phi (\var)[\var'] = \p'(T)$. Therefore, since
$J(\varepsilon)=\langle\p(T)|O|\p(T)\rangle-\alpha\ds\int_0^T\varepsilon^2(t)dt$,
we have  
$$(\nabla J (\var),\var')= 2 ~Re \langle\p'(T)|O|\psi(T)\rangle -2 \alpha \int_0^T  \var (t) \var' (t)dt.$$
To end the proof of (\ref{dJ}), we consider the solution $\chi$ of the adjoint
state equation (\ref{schrodadj}) and we multiply equation (\ref{eqschrod1}) by
$\overline\chi$ (the complex conjugate of $\chi$), integrate on $\mathbb R
\times [0,T]$ and take the imaginary part. We obtain: 
$$Im \int_0^T \int_{\mathbb{R}} (i\partial_t \p' -[H-\mu\var]\p' )\overline{\chi}
=Im \int_0^T \int_{\mathbb{R}}\mu \var'\p \overline{\chi}.$$
After an integration by parts and since $\p'(0)=0$, we get
$$Im \int_0^T \int_{\mathbb{R}} ~ \overline{i\partial_t\chi}  \p'- Im
\int_{\mathbb{R}} \p'(T)~ 
\overline{i\chi(T)} - Im \int_0^T \int_{\mathbb{R}}~ \overline{[H-\mu\var]\chi} \p' 
=Im \int_0^T \int_{\mathbb{R}}\mu \var'\p \overline{\chi}.$$
Since $\chi$ satisfies equation~(\ref{schrodadj}), we then obtain
$$Re \langle\p'(T)|O|\psi(T)\rangle = Re  \int_{\mathbb{R}} \p'(T)~  \overline{O\p(T)}  
=- Im \int_0^T\!\! \int_{\mathbb{R}}\mu \var'\p \overline{\chi}
= -\int_0^T \!\!Im\langle\chi(t)|\mu|\psi(t)\rangle \var' (t)dt$$ 
what ends the proof of the lemma.

\qed

In what follows, we denote by $\nabla J (\var)$ the function
$~t \mapsto -2\big(\alpha \var (t) + Im\langle\chi(t)|\mu|\psi(t)\rangle\big)~$ 
and by $~C_J~$ the set of the critical points of $J$, \textit{i.e.},
\beq\label{defCJ}
C_J=\big\{\var\in L^2(0,T), \quad \forall t\in [0,T],\quad\alpha\var (t) +Im \langle \chi (t) |\mu|\psi(t)\rangle=0 \big\}.
\eeq
Note that, thanks to the results of the section \ref{sec2}, we have $C_J\subset
L^\infty (0,T)$ since for all $\var\in C_J$, 
$$\|\var\|_{L^\infty(0,T)} ~\leq ~\dfrac{1}{\alpha} 
\|\langle \chi |\mu|\psi\rangle\|_{L^\infty(0,T)}
~\leq~ C \|\mu\|_{L^\infty}\|\chi\|_{L^\infty(0,T;L^2)}\|\psi\|_{L^\infty(0,T;L^2)}.$$

\begin{remark}\label{rem:singlet}{\bf :} 
Note also that for $\alpha>6 T\|\mu\|^2_{L^\infty}\|O\|_*$, the set $C_J$ is reduced to one point. Indeed, suppose
that $C_J$ contains two distinct points $\var_1$ and $\var_2$, we then have, for
$t\in (0,T)$:
\beq\nonumber
\alpha\big(\var_2(t)-\var_1(t)\big)+Im\langle
\chi_2(t)-\chi_1(t)|\mu|\psi_2(t)\rangle+Im\langle\chi_1(t)|\mu|\psi_2(t)-\psi_1(t)\rangle =0,
\eeq
 where $\psi_1$, $\psi_2$ (resp. $\chi_1$, $\chi_2$) are the solutions of
 $(\ref{schrod0eq})$ (resp. $(\ref{schrodadj})$)
 corresponding to $\var_1$ and $\var_2$ respectively.
Using estimates $(\ref{estimp})$ with $\psi=\psi_1$, $\psi'=\psi_2-\psi_1$,
$\var=\var_2$ and $\var'=\var_2-\var_1$ and $(\ref{estimadj})$ with
$\chi=\chi_1$, $\chi'=\chi_2-\chi_1$, $\psi'=\psi_2-\psi_1$, $\var=\var_2$ and
$\var'=\var_2-\var_1$, we obtain
\beq\nonumber
\alpha\|\var_2-\var_1\|_{L_1(0,T)}\leq 6 T\|\mu\|^2_{L^\infty}\|O\|_*\|\var_2-\var_1\|_{L^1(0,T)},
\eeq
which leads to
$\alpha\leq 6 T\|O\|_*\|\mu\|^2_{L^\infty}$, 
 and the result follows.
\end{remark}

In order to prove the compactness of $C_J$, we introduce an important
property of the application $\var(t) \mapsto \p(x,t)$, firstly
presented in a more general setting by J. M. Ball, J. E. Marsden and
M. Slemrod in \cite{B-M-S}. In our context, this result can be stated as
follows.  

\begin{lemma}\label{ballslemrod}
Assume that $\var\in L^1(0,T)$, $\mu : X\to X$ is a bounded operator and that
$H$ generates a $C^0$-semigroup of bounded linear operators on some Banach
space $X$. For $x\in \bbR^3$ and $t \in (0,T)$, we denote by $\p(x,t)$ the solution of  
\beq\nonumber
\left\{\begin{array}{l}
i\partial_t \psi -[H-\mu\var]\psi =0,\\
\psi(0)=\p_0\in X.
\end{array}
\right.
\eeq
Then, $\var \mapsto \p$ is a compact mapping in the sense that for any weakly converging sequence $(\var_n)_{n\in\bbN}$ to $\var$ in $L^1(0,T)$, $(\p_n)_{n\in\bbN}$ converges strongly to $\p$ in $C([0,T];X)$.
\end{lemma}
The precise proof of this result derives directly from \cite{B-M-S} (Theorem 3.6, p580), see also  \cite{bookgt}  and  \cite{these}. 
It allows us to obtain the following lemma.

\begin{lemma}\label{compacite}
For $\mu\in W^{2,\infty}$, 
$C_J$ is compact in $L^\infty(0,T)$.
\end{lemma}
{\bf Proof:} \indent
Consider a bounded sequence $(\var^n)_{n\in\bbN}$ of $C_J$. 
By definition, for all $n\in\bbN$, $\var^n\in L^2(0,T)$ and  $~\var^n(t) =
-\dfrac{1}{\alpha}~ \langle \chi^n(t) |\mu|\psi^n(t)\rangle~$ where $\psi^n$ and 
$\chi^n$ are the corresponding solutions of $(\ref{schrod0eq})$ and $(\ref{schrodadj})$.
It is also possible to extract  a weakly convergent sub-sequence in $L^2(0,T)$, still denoted $(\var^n)_{n\in\bbN}$.
From Lemma \ref{schrod0}, one knows that the Hamiltonian 
$H=-\Delta + V$ with $V\in W^{2,\infty}$ generates a $C^0$-semigroup of bounded linear operators on the Banach space $X = H^2$. Therefore, with  $\mu\in W^{2,\infty}$ the conditions of Lemma~\ref{ballslemrod} are fulfilled and we obtain the strong convergences $\psi^n\stackrel{n\to +\infty}{\longrightarrow} \p$ and 
$\chi^n\stackrel{n\to +\infty}{\longrightarrow} \chi$ in $C([0,T];H^2)$.
Thus, for all $t\in(0,T)$,
$\displaystyle\int_{\bbR^3} \overline{\psi^n(t)}~\mu\chi^n(t)\, dx \stackrel{n\to +\infty}{\longrightarrow} \displaystyle\int_{\bbR^3} \overline{\psi(t)}~\mu\chi(t)\, dx.$
The sequence $\left(\var_n(t)\right)_{n\in\bbN}$ then strongly converges in $L^\infty(0,T)$ and the result follows.

\qed
\subsection{Analyticity of $J$}
The implicit formulation of the derivative can be iteratively carried on in order to prove
the analyticity of $J$.

\begin{lemma}\label{analiticity}
Let $\psi$ be the solution of $(\ref{schrod0eq})$ corresponding to $\varepsilon$.
The functional 
\begin{eqnarray*}
 \vartheta : L^2(0,T)&\to& L^\infty(0,T;H^2)\cap W^{1,\infty}(0,T;L^2)\\
\var&\mapsto&\psi,
\end{eqnarray*}
 is analytic.
\end{lemma}
{\bf Proof:} \indent Let $\varepsilon,\varepsilon' \in
L^2(0,T)$ be such that $\|\varepsilon'\|_{L^1(0,T)}\leq
\dfrac{1}{4\|\mu\|_{L^\infty}}$ and the sequence
$(\p^\ell)_{\ell\in \bbN}\in (L^\infty(0,T;H^2))^{\bbN}$ defined recursively by
$\p^0=\vartheta(\varepsilon)$ and for $\ell>0$: \beq\label{schroddern}
\left\{\begin{array}{l}
 i\partial_t  \psi^\ell(x,t) -[H - \mu(x) \varepsilon(t)]\psi^\ell(x,t) =-  \mu(x) \varepsilon'(t)\psi^{\ell-1}(x,t) \\
 \psi^{\ell}(x,0) =0 .
\end{array}\right.
\end{equation}
The existence of $\psi^\ell$ is consequence of Lemma~\ref{schrod1}. Thanks to
(\ref{estimp}) applied with $\psi=\psi^{\ell-1}$ and $\psi'=\psi^\ell$, one has for $\ell\geq 1$ and $t\in [0,T]$:
\beq\label{gron1}
\|\p^\ell(t)\|_{L^2}\leq2\|\mu\|_{L^\infty}\|\var'\|_{L^1(0,T)}\|\p^{\ell-1}\|_{L^\infty(0,T;L^2)}
\leq 2^\ell \|\mu\|_{L^\infty}^{\ell}
\|\varepsilon'\|^\ell_{L^1(0,T)}. \eeq 
Given $N>0$, we obtain by summing (\ref{schroddern}) from $\ell=0$ to $N$:
\beq\label{schroddernsum}
\left\{\begin{array}{l}
 i\partial_t  \left(\ds\sum_{\ell=0}^N\psi^\ell(x,t) \right)-
 [H - \mu(x) \big(\varepsilon(t)+\varepsilon'(t)\big)]
 \left(\ds\sum_{\ell=0}^N\psi^\ell(x,t) \right) =\mu(x) \varepsilon'(t)\psi^{N}(x,t) \\ 
 \ds\sum_{\ell=0}^N\psi^\ell(x,0)=\p_0(x).
\end{array}\right.
\end{equation}
On the other hand, one has:
\beq\label{schroddl}
i\partial_t\vartheta(\var+\var')-[H - \mu(x) (\varepsilon(t)+\varepsilon'(t))]\vartheta(\var+\var')=0.
\eeq
Subtracting (\ref{schroddernsum}) and (\ref{schroddl}) and using estimates
(\ref{estimp}) with $\psi=-\psi^{N}$,
$\psi'=\sum_{\ell=0}^N\psi^\ell(x,t)-\vartheta(\var+\var')$, $\var=\var+\var'$
and (\ref{gron1}), we get: 
\beq\nonumber
\left\|\vartheta(\varepsilon+\varepsilon')(t)-\sum_{\ell=0}^{N}\psi^\ell(t)\right\|_{L^2}\leq 2^N
\|\mu\|_{L^\infty}^{N} \|\varepsilon'\|^N_{L^1(0,T)} \leq
2^{-N}\eeq
and the functional $\vartheta$ reads now: $
\vartheta(\varepsilon+\varepsilon')=\ds\sum_{\ell=0}^{\infty}\psi^\ell$ in $L^2(0,T)$.
Since
$\varepsilon'\mapsto \psi^{\ell}$ is $\ell$-linear, the theorem follows. 

\qed
The next lemma follows immediately from this result.

\begin{lemma}
The cost functional $~J$ is analytic
on $L^2(0,T)$. 
\end{lemma}

\subsection{About the Hessian operator of $J$}
Let us now investigate some properties of the second order derivative of
$J$. Though we express it as an implicit function of its argument $\var$, some
results can be obtained from the next lemma.

\begin{lemma}\label{hes}
Let $\p$ and $ \chi$ be the solutions of $(\ref{schrod0eq})$ and
$(\ref{schrodadj})$.
The functional $\gamma : \varepsilon \mapsto  Im
\langle\chi|\mu|\psi\rangle$ is differentiable on
$L^2(0,T)$ and one has:
\beq
\begin{array}{cccl}\label{deriv}
D\gamma(\varepsilon)[\var']= Im
\langle\chi'|\mu|\psi\rangle+Im\langle\chi|\mu|\psi'\rangle,
\end{array}
\eeq
where $\psi'$ and $\chi'$ are the solutions of
$(\ref{eqschrod1})$ and $(\ref{schrodadj1})$. Moreover, for all $\var\in
L^2(0,T)$, $D\gamma(\var)$ is compact on $L^2(0,T)$. 
\end{lemma}
{\bf Proof:} \indent
Let $\var\in L^2(0,T)$ and $\psi$ and $\chi$ the corresponding solutions of
(\ref{schrod0eq}) and (\ref{schrodadj}).
As in the proof of Lemma~\ref{grad}, the key-point is the differentiability of 
the functional $\vartheta$, defined in Lemma \ref{analiticity} on $L^2(0,T)$.  
Actually, $D\vartheta(\var)[\var'] = \p'$, where $\p'$ is the solution of (\ref{eqschrod1}). 
The main explanations can be read 
in \cite{baudouin2}. Repeating this argument for $\var\mapsto\chi$, we obtain that $\gamma$ is differentiable and we get
(\ref{deriv}).\\ 
\indent Let us now prove the compactness of this operator. Let
$(\var'^n)_{n\in\bbN}$ be a bounded sequence in $L^2(0,T)$ and let
$(\psi'^n)_{n\in\bbN}$ and  $(\chi'^n)_{n\in\bbN}$  be the
corresponding solutions of (\ref{eqschrod1}) and
(\ref{schrodadj1}). \\
As $\psi'^n\in L^\infty(0,T;H^2)\cap W^{1,\infty}(0,T;L^2) $  (see the proof of
Lemma~\ref{schrod1}), we have that $\psi'^n \in C([0,T];L^2)$ and $\partial_t\psi'^n \in
L^2(0,T;L^2)$. By means of the continuity of:
\beq
\begin{array}{ll ll lll }
L^2(0,T)&\rightarrow& C([0,T];L^2)&\textnormal{and}&L^2(0,T)&\rightarrow& L^2(0,T;L^2)\\
\hfill{}\var'&\mapsto& \psi'&&\hfill{}\var'&\mapsto& \partial_t\psi'
\end{array}
\eeq
there exist $\psi'^\infty$ such that, up to extraction,
$ \psi'^n \rightharpoonup \psi'^\infty \in L^2(0,T;L^2)$ and
\beq \partial_t\psi'^n \rightharpoonup \partial_t \psi'^\infty \in L^2(0,T;L^2).\label{estim2}
\eeq
Since $\psi'^n(0)=0$, we have
$\psi'^n(t)=\displaystyle\int_0^t \partial_t\psi'^n(s)ds$ and \eqref{estim2} implies that for all $t\in [0,T]$, $(\|\psi'^n(t)\|_{L^2})_{n\in\bbN}$ is uniformly bounded. 
Moreover, for all $t,t'\in [0,T],
t\leq t'$, we have:
$$ \|\psi'^n(t')-\psi'^n(t)\|_{L^2}\leq
\int_t^{t'}\|\partial_t\psi'^n(s)\|_{L^2}ds
\leq\sqrt{t'-t}\|\partial_t\psi'^n\|_{L^2(0,T;L^2)}.$$
Combining this with \eqref{estim2}, we find that
$(\psi'^n)_{n\in\bbN}$ is an equicontinuous sequence in $C([0,T],L^2)$. 
We conclude by applying Ascoli's theorem to the family $\left\{ Im \langle \chi|\mu|\psi'^n\rangle,n\in \bbN \right\}$ of  the space $C([0,T])$. 
Similar arguments apply for
$\left\{Im\langle\chi'^n|\mu|\psi\rangle,n\in \bbN \right\}$, and the
results follows. 

\qed

Thanks to the previous lemma, $J$ is twice differentiable and its Hessian
operator reads:
\beq\nonumber
H_J(\varepsilon):\varepsilon'\mapsto -2\big(\alpha \varepsilon' +
D\gamma(\varepsilon)[\varepsilon']\big).
\eeq
In the sequel, a criterion ensuring that the Hessian operator of $J$ is invertible will be
useful. The next lemma provides it.

\begin{lemma}\label{esthess}
Suppose that:
$\alpha>6T\|\mu\|^2_{L^\infty}\|O\|_*.$
Then the operator $H_J(\var)$ is invertible on $L^2(0,T)$.
\end{lemma}
{\bf Proof:}\indent We keep the notation of Lemma \ref{hes}. The Cauchy-Schwarz inequality, combined with (\ref{estimp}) and (\ref{estimadj2}) yields:
\beqn\nonumber
\|D\gamma(\varepsilon)[\varepsilon']\|_{L^\infty(0,T)}&\leq&
\|\mu\|_{L^\infty}\big(\|\chi'\|_{L^\infty(0,T;L^2)}+\|O\|_*\|\psi'\|_{L^\infty(0,T;L^2)}\big)\\\nonumber &\leq&
6\sqrt T\|\mu\|^2_{L^\infty}\|O\|_*\|\var'\|_{L^2(0,T)}.
\eeqn
Finally, thanks to the assumption of the lemma, one has
\beq\nonumber
\sup_{\{\varepsilon',
  \|\varepsilon'\|_{L^2(0,T)}=1\}}\left(\frac{1}{\alpha}\|D\gamma(\varepsilon)[\varepsilon']\|_{L^2(0,T)}\right)< 1, 
\eeq
which implies that $I+\dfrac 1\alpha D\gamma(\varepsilon)$ is invertible and the result follows.

\qed

\section{\L ojasiewicz inequality for the cost functional $J$}\label{sec5}
Several convergence results of dynamical systems have been proved thanks to the
\L ojasiewicz inequality recalled here. In order to tackle the problem
of the convergence of the optimizing sequence presented in the next
section, we have to extend this 
inequality to the case of a compact set in an infinite dimensional space.
The basic result considered in this section is the following
(cf \cite{loja,loja2}):

\begin{theorem}\label{loja}
Let $N$ be an integer and $\Gamma:\bbR^N\rightarrow \bbR$ be an
analytic function in a neighborhood of a point $a\in\bbR^N$. Then
there exists $\sigma>0$ and $\theta\in ]0,\frac{1}{2}]$ such that
\beq\label{eqloja}
\forall x\in \bbR^N, \ \|x-a\|<\sigma,\qquad \|\nabla \Gamma
(x)\|\geq |\Gamma(x)-\Gamma(a)|^{1-\theta}, 
\eeq
where $\|.\|$ is a given norm on $\bbR^N $.
\end{theorem}
The real number $\theta$ is a \L ojasiewicz exponent of $a$.
Following the work \cite{jendoubi} of M. A. Jendoubi (which
simplifies the theorem of \L ojasiewicz-Simon \cite{simon}), the
latter theorem can be generalized to the case of infinite dimension.

\begin{lemma}\label{loja2}
Given $\var\in L^2(0,T)$, there exists $\sigma'>0$, $\kappa>0$ and $\theta' \in ]0,\frac{1}{2}]$
  such that:
\beq\nonumber
\forall \var'\in  L^2(0,T),
  \|\var'-\var\|_{L^2(0,T)}\leq \sigma' ,\quad \|\nabla
  J(\var')\|_{L^2(0,T)}\geq \kappa
 |J(\var')-J(\var)|^{1-\theta'}.
\eeq
\end{lemma}
We give the proof of this lemma in the appendix.
A more precise result can be obtained if the Hessian operator under
consideration is invertible at point $a$ (see e.g, \cite{kavian}).
Indeed, one can then show that $1/2$ is a \L ojasiewicz exponent of $a$. 
We will use this improvement in Section~\ref{sec8} since  Lemma~\ref{esthess} provides actually an expected sufficient condition. 
The next lemma is a global version of the previous one.

\begin{lemma}\label{ttk}
Let $\widetilde{C}_J$ be a connected component of
$C_J$ in $L^2(0,T)$. We denote by $l$ the value of $J(\var)$ for all $\var\in \widetilde{C}_J$ and we set $\widetilde{J}(\varepsilon)=l-J(\varepsilon)$. 
There exist $\widetilde{\sigma}>0$, $\wk>0$ and
$\wt\in ]0,\frac{1}{2}]$ such that:
 \beq\label{estloja} \forall \varepsilon \in L^2(0,T),\
d_2(\varepsilon,\widetilde{C}_J)<\widetilde{\sigma}, \qquad \|\nabla
J(\varepsilon)\|_{L^2(0,T)}\geq \wk
|\widetilde{J}(\varepsilon)|^{1-\wt} ,\eeq where $d_2$ is the
distance associated to the $L^2(0,T)$-norm.
\end{lemma}
{\bf Proof:}\indent
Lemma~\ref{loja2} ensures that for each point $a$ in $\widetilde{C}_J$
there exist three real numbers  $\sigma_a$, $\theta_a$ and $\kappa_a$ such that:
\beq\nonumber \forall \var\in \bbR^N,\ \|\var-a\|_{L^2(0,T)}<\sigma_a \qquad \|\nabla
J (\var)\|_{L^2(0,T)} \geq \kappa_a|\widetilde{J}(\var)|^{1-\theta_a}.\eeq The
compactness of $\widetilde{C}_J$, guaranteed by Lemma~\ref{compacite}, 
allows us to extract from $\left\{B(a,\frac{\sigma_a}{2}), a\in \widetilde{C}_J\right\}$
a finite family $A=\left\{B(a_i,\frac{\sigma_{a_i}}{2})\right\}_{i\in
F}$, where $F$ is a finite set of indexes, such that $\widetilde{C}_J\subset A$. 

We then define $\ws$, $\wk$ and $\wt \in ]0,1/2]$ as the respective lower
bounds of $\left\{\frac{\sigma_{a_i}}{2}\right\}_{i\in F}$,
$\left\{\kappa_{a_i}\right\}_{i\in F}$ and $\left\{\theta_{a_i}\right\}_{i\in 
F}$ and the result follows.

\qed

\section{Optimizing sequence}\label{sec6}
We have now gathered all the necessary results to present and analyze the optimizing sequence.  

\subsection{Definition of the sequence}
Following the approach sketched in the introduction, Y. Maday
and G. Turinici have defined an optimizing
sequence $(\varepsilon^k)_{k\in\bbN}$ for the cost functional $J$ as follows \cite{maday} :\\

Consider $(\delta,\eta)\in ]0,2[\times]0,2[$, $\varepsilon^0\in
    L^{\infty}(0,T)$, $\et^0\in
    L^{\infty}(0,T)$, $\p^0$ and $\chi^0$ the corresponding solutions of (\ref{schrod0eq})
    and (\ref{schrodadj}) according to Lemma~\ref{schrod0}. The functions $\varepsilon^{k}$ and $\et^{k}$ are computed by
    solving iteratively:
\begin{eqnarray} \label{algo1}
& \! & \left\{
\begin{array}{l}
  i \partial_t \p^{k}(x,t) =
 \big( H (x)-\varepsilon^k(t)\mu (x)\big) \p^{k}(x,t)
\\
 \psi^k (x,0) = \psi_0(x)
\end{array} \right.
\\ \label{algo2}
& \ & \varepsilon^k(t) = (1-\delta)\tilde{\varepsilon}^{k-1}(t) -
\frac{\delta} {\alpha} Im \langle \chi^{k-1} (t)| \mu | \p^{k}(t)
\rangle
 \\ \label{algo3}
& \! & \left\{
\begin{array}{l}
  i \partial_t \chi^{k}(x,t) =
\big(H (x)-  \tilde{\varepsilon}^k(t)\mu (x)\big)%
 \chi^k(x,t)  \\
 \chi^k (x,T) = O \p^k(x,T)
\end{array} \right.
\\ \label{algo4}
&  \ & \tilde{\varepsilon}^k(t) =  (1-\eta)\varepsilon^{k}(t) -
\frac{\eta}{\alpha} Im   \langle \chi^{k}(t) | \mu | \p^{k}(t)
\rangle.
\end{eqnarray}
Existence and uniqueness of solutions $\p^k$ and $\chi^k$ of the above
equations result from an easy adaptation of Lemma~\ref{schrod2}, as for the proof of $\var^k, \tilde{\varepsilon}^k\in L^2(0,T)$ for all $k\in\bbN$.\\
\begin{remark}{\bf :}
Note that this choice of optimizing sequence is not canonical. There
exists other ways to guarantee that the condition $(\ref{mono-cond})$
is fulfilled (see, {\it e.g.} $\cite{tur_cdc}$). However, this
formulation includes many monotonic algorithms, {\it e.g.} the one by
Krotov (presented in $\cite{tannor}$) or by W. Zhu and H. Rabitz
$\cite{zhu}$ which are often used in the numerical simulations. 
\end{remark}

\subsection{Properties of the sequence}
We present here two results about the sequence $(\var^k)_{k \in
  \bbN}$. The proofs can be found in \cite{salomon,maday}. These results
state that $(\varepsilon^k)_{k\in\bbN}$ defined by
$(\ref{algo1})-(\ref{algo4})$ is bounded in $L^\infty(0,T)$ and that the
corresponding sequence $\big(J(\var^k)\big)_{k\in\bbN}$ increases
monotonically.
\begin{lemma}\label{lem:borne}
Given an initial field $\varepsilon^0\in L^\infty(0,T)$, 
let us define $M$ by:
\begin{equation}\nonumber
M=\max\left(\|\varepsilon^0\|_{L^\infty(0,T)},\max\left(1,\frac{\delta}{2-\delta},\frac{\eta}{2-\eta} \right)\frac{\|O\|_{*}\|\mu\|_{L^\infty}}{\alpha}\right).
\end{equation}
The sequences $(\varepsilon^k)_{k\in\bbN}$ and
$(\widetilde{\var}^k)_{k\in\bbN}$ satisfy: 
\beq\nonumber
\forall k\in
\bbN,\ \|\varepsilon^k\|_{L^\infty(0,T)}\leq M,\
\|\et^k\|_{L^\infty(0,T)}\leq M. \eeq
\end{lemma}
\begin{lemma}\label{lem:monotonic}
The sequence $(\varepsilon^k)_{k\in\bbN}$ defined by
$(\ref{algo1})--(\ref{algo4})$ ensures the monotonic convergence of the
cost functional $J$ in the sense that: \beqn  \nonumber
J(\varepsilon^{k+1})-J(\varepsilon^{k})&=&\langle\p^{k+1}(T)-\p^{k}(T)|O|\p^{k+1}(T)-\p^{k}(T)\rangle\\
 &&+\left(\frac{2}{\eta}-1\right)\left\|\varepsilon^{k+1}-\et^k\right
 \|^2_{L^2(0,T)}+\left(\frac{2}{\delta}-1\right)\left\|\et^{k}-\varepsilon^k\right\|^2_{L^2(0,T)}, \label{monotonic}
\eeqn and there exists $l_{\varepsilon^0}$ such that $
\ds\lim_{k\rightarrow +\infty} J(\varepsilon^{k})=l_{\varepsilon^0}$.
\end{lemma}
In order to study the convergence of $(\varepsilon^k)_{k\in\bbN}$, we will need to estimate the gradient of  $J$ at each point $\varepsilon^k$. Such an estimate is obtained in the next lemma.

\begin{lemma}\label{estim}
There exists $\lambda>0$, depending on $\mu$, $O$, $\alpha$,
$\delta$, $\eta$ and $T$, such that:
\beq\label{estnab}
\|\nabla J(\varepsilon^k)\|_{L^1(0,T)}\leq\lambda
\left(\|\varepsilon^{k}-\et^{k-1}\|_{L^2(0,T)}+\|\et^{k-1}-\varepsilon^{k-1}\|_{L^2(0,T)}\right).
\eeq
\end{lemma}
{\bf Proof:}\indent
Thanks to (\ref{dJ}), we have:
\beqn
\nabla J(\varepsilon^k)(t)&=&-2\big(\alpha \var^k (t) + Im
\langle\chi^{k-1}(t)|\mu|\psi^k(t)\rangle+Im
\langle\chi^{\var^k}(t)-\chi^{k-1}(t)|\mu|\psi^k(t)\rangle\big)\nonumber\\&&\nonumber\\
&=&-2\left(\alpha \left(1-\dfrac{1}{\delta}\right)\big( \var^k(t)-\widetilde{\var}^{k-1} (t)\big)+Im
\langle\chi^{\var^k}(t)-\chi^{k-1}(t)|\mu|\psi^k(t)\rangle \right)\label{der1}
\eeqn
where $\chi^{\var^k}$ is the solution of (\ref{schrodadj}) with $\psi=\psi^k$
and $\var=\var^k$.\\ 
Next, $\chi^{\var^k}-\chi^{k-1}$ is the solution of equation (\ref{schrodadj1}) corresponding to  $\var=\var^k$, $\var'=\widetilde{\var}^{k-1}-\var^k$
and $\chi=\chi^{k-1}$. The associated estimate (\ref{estimadj}) then
gives:
\beqn
\|\chi^{\var^k}(t)-\chi^{k-1}(t)\|_{L^2}
&\leq& 2 \|\mu\|_{L^\infty}\|\var^{k}-\widetilde{\var}^{k-1}\|_{L^1(0,T)}\|O\|_*
+\|O(\p^{k}(T)-\p^{k-1}(T))\|_{L^2}\nonumber\\&&\nonumber\\
&\leq& 4\|\mu\|_{L^\infty}\|O\|_*(\|\var^k-\widetilde{\var}^{k-1}\|_{L^1(0,T)}+
\|\widetilde{\var}^{k-1}-\var^{k-1}\|_{L^1(0,T)}).\label{der2}    
\eeqn
Combining (\ref{der1}) and (\ref{der2}), we obtain (\ref{estnab}) with
$\lambda=2\sqrt{T}\left(4T\|O\|_*\|\mu\|^2_{L^\infty}+\alpha \left(1-\dfrac{1}{\delta}\right)\right)$.

\qed

\subsection{Limit points of the sequence}
We now present some result about the limit points of 
$(\varepsilon^k)_{k\in\bbN}$. These results give first hints about the
relationship between these limit points and the set $C_J$ of the critical points of
the cost functional $J$. Thus, we obtain a first case of convergence.
\begin{lemma}
Let $(\varepsilon^{k_n})_{n\in\bbN}$ be a weakly convergent
sub-sequence of $(\varepsilon^{k})_{k\in\bbN}$ in $L^2(0,T)$.
Then $(\varepsilon^{k_n})_{n\in\bbN}$ converges in
$L^\infty(0,T)$ towards a critical point of the cost functional $J$.
\end{lemma}
{\bf Proof:} \indent
Let $(\varepsilon^{k_n})_{n\in\bbN}$ be a weakly convergent
sub-sequence of $(\varepsilon^{k})_{k\in\bbN}$ in $L^2(0,T)$ and let us consider $\ell \in \bbN$. Equation
(\ref{monotonic}) ensures that $(\varepsilon^{k_n+\ell})_{n\in\bbN}$ also converges
weakly (and has the same limit as $(\varepsilon^{k_n})_{n\in\bbN}$).
Thanks to Lemma~\ref{ballslemrod}, the sequences
$\left(\chi^{k_n+\ell}\right)_{n\in\bbN}$ and $\left(\p^{k_n+\ell}\right)_{n\in\bbN}$ converge strongly in the space
$C([0,T];L^2)$.  Thus, we obtain by bilinearity the strong
convergence of both sequences $\left(\left\langle \chi^{k_n} |\mu
|\psi^{k_n} \right\rangle\right)_{n\in\bbN}$ and $\left(\left\langle
\chi^{k_n} |\mu |\psi^{k_n+1} \right\rangle \right)_{n\in\bbN}$ in  
$L^\infty (0,T)$.\\

According to (\ref{algo2}) and (\ref{algo4}), $(\var^{k_n})_{n\in\bbN}$ also reads:
\beq\nonumber
\var^{k_n+1}=\underbrace{(1-\delta)(1-\eta)}_{\nu}\var^{k_n}+u_{k_n},
\eeq
where $|\nu|<1$ and where
$~ u_{k_n}(t)=-\dfrac {(1-\delta)\eta}\alpha Im \left\langle
\chi^{k_n}(t)|\m|\p^{k_n}(t)\right\rangle-\dfrac {\delta}\alpha Im \left\langle
\chi^{k_n}(t)|\m|\p^{k_n+1}(t)\right\rangle~$
strongly converges in $L^\infty (0,T)$. Note again that given $\ell\in \bbN$,
$(u_{k_n+\ell})_{n\in \bbN}$ also converges in $L^\infty(0,T)$ (towards the
same limit).
For all $k\in\bbN$, the absolute value of $u_k(t)$
can be estimated by:
\beq\nonumber
|u_k(t)|\leq m=\dfrac{4\|\mu\|_{L^\infty}\|O\|_*}\alpha.
\eeq
Let us prove that $(\varepsilon^{k_n})_{n\in\bbN}$ is Cauchy in $L^\infty(0,T)$.
Consider $e>0$.
There exists $n_1>0$ be such that
\begin{equation}\label{cn'}
2m\sum_{j_1}^\infty |\nu|^j \leq \frac{e}{4}.
\end{equation}
Since the sequence $(u_{k_n-\ell})_{n\in\bbN}$ is Cauchy for all $\ell$ with
$0\leq\ell\leq n_1$, we have:
\begin{equation}\label{Cauchy}
\exists n_2 >0/ \ \forall s > n_2, \ \forall q\geq 0, \
\|u_{k_{s+q}-\ell}-u_{k_s-\ell}\|_{L^\infty(0,T)}\leq \frac{e}{4n_1}.
\end{equation}
Let $n$ be an integer fulfilling the conditions:
\begin{equation}\label{cn1}
\forall p\geq 0,\ |\nu^{k_{n+p}}-\nu^{k_{n}}|\leq\frac{e}{4\|\var^0\|_{L^\infty(0,T)}},\quad  k_n>n_1,\quad n>n_2.
\end{equation}
Let $p$ be a positive integer. Since we have, for all $n\in\bbN^*$,
$\var^{k_n}=\nu^{k_n}\var^{0}+\ds\sum_{j=0}^{k_n-1}\nu^j u_{k_n-j-1}$
we obtain
\begin{equation}\label{cauch}
\begin{array}{ccl}
\var^{k_{n+p}}-\var^{k_n}&=&(\nu^{k_{n+p}}-\nu^{k_n})\var^0+
\ds\sum_{j=k_n}^{k_{n+p}-1}\nu^{j}u_{k_{n+p}-j-1}\\
&&+\ds\sum_{j_1}^{k_n-1}\nu^{j}(u_{k_{n+p}-j-1}-u_{k_n-j-1})
+\ds\sum_{j=0}^{n_1-1}\nu^{j}(u_{k_{n+p}-j-1}-u_{k_n-j-1}).
\end{array}
\end{equation}
Thank to (\ref{cn'}) and the two first conditions of (\ref{cn1}):
\beq\nonumber
\left\|(\nu^{k_{n+p}}-\nu^{k_n})\var^0\right\|_{L^\infty(0,T)}\leq\frac{e}{4},
\eeq
\beq\nonumber
\left\|\sum_{j=k_n}^{k_{n+p}-1}\nu^{j}u_{k_{n+p}-j-1}\right\|_{L^\infty(0,T)}
\leq\sum_{j=k_n}^{\infty}\left\|\nu^{j}u_{k_{n+p}-j-1}\right\|_{L^\infty(0,T)}\leq
m\sum_{j=k_n}^{\infty}|\nu|^{j}\leq\frac{e}{4}.
\eeq
According to the condition (\ref{cn'}), the third term of (\ref{cauch})
can be estimated by:
\beq\nonumber
\left\|\sum_{j_1}^{k_n-1}\nu^{j}(u_{k_{n+p}-j-1}-u_{k_n-j-1})\right\|_{L^\infty(0,T)}
\leq 2m\sum_{j_1}^{\infty}|\nu|^{j}\leq\frac{e}{4}.
\eeq
Lastly, $|\nu|<1$, the third condition of (\ref{cn1}) and the Cauchy property (\ref{Cauchy}) 
 allows us to estimate the last term of (\ref{cauch}):
\begin{eqnarray*}
\left\|\sum_{j=0}^{n_1-1}\nu^{j}(u_{k_{n+p}-j-1}-u_{k_n-j-1})\right\|_{L^\infty(0,T)}
&\leq& \sum^{n_1-1}_{j=0}\left\|u_{k_{n+p}-j-1}-u_{k_n-j-1}\right\|_{L^\infty(0,T)}~\leq~\frac{e}{4}.
\end{eqnarray*}
We have thus proved that for all $e>0$, if $n$ is large enough then, for every $p>0$,
\beq\nonumber
\|\var^{k_{n+p}}-\var^{k_{n}}\|_{L^\infty(0,T)}\leq e,
\eeq
which proves that $(\var^{k_n})_{n\in\bbN}$ is Cauchy in $L^\infty(0,T)$.\\
\indent We denote by $\var$ the limit of $(\var^{k_n})_{n\in\bbN}$. Thanks to
(\ref{monotonic}), $(\et^{k_n})_{n\in\bbN}$ also converges towards $\var$.
Passing through the limit in $(\ref{algo1})-(\ref{algo4})$,
we then deduce that $\var$ belongs to $C_J$, according to definition~(\ref{defCJ}).

\qed

Let us denote by $C_{\varepsilon^0}\subset C_J$ the set of the
limit points of $(\varepsilon^{k_n})_{n\in\bbN}$. As stated in Remark
\ref{rem:singlet}, for  $\alpha>6 T\|\mu\|^2_{L^\infty}\|O\|_*$, $C_J$, and
consequently $C_{\varepsilon^0}$, are
reduced to one point. By means of Lemma \ref{lem:borne}, the convergence of the sequence
$(\varepsilon^k)_{k\in\bbN}$ is then guaranteed in this case. 
\begin{remark}\label{rem:conv_mini}{\bf :}
In addition, the uniqueness of the critical point implies that   
the limit in this case is necessarily an extremum of $J$.
\end{remark}
In order to obtain the convergence for all $\alpha>0$, we need to study more
precisely the asymptotic behavior of the sequence
$(\varepsilon^k)_{k\in\bbN}$ in the neighborhood of $C_{\varepsilon^0}$.
A standard argument of compactness applied to $C_{\varepsilon^0}$
enables us to obtain the following result.

\begin{lemma}
Let denote by $d_\infty$ the distance corresponding to the
$L^\infty(0,T)$ norm. One has: \beq\label{dist}
d_\inftyÒ(\varepsilon^k,C_{\varepsilon^0}) \rightarrow 0. \eeq
\end{lemma}

\begin{remark}{\bf :}
By means of the monotonicity property, we find that $J=l_{\varepsilon^0}$ on the set $C_{\varepsilon^0}$ (with $l_{\varepsilon^0} = \lim_{k\rightarrow +\infty} J(\varepsilon^{k})$). It is then possible to apply Lemma~\ref{ttk} with
$\widetilde{C}_J=C_{\var^0}$ since the assumption that $\widetilde{C}_J$ is connected is
only necessary to ensure that $J$ is constant on this set. It can however be proved that
$C_{\varepsilon^0}$ is connected (see \cite{salomon-cdc05}).
\end{remark}

\section{Convergence of the sequence}\label{sec7}
It is now possible to prove the convergence of the sequence
$(\varepsilon^k)_{k\in\bbN}$ by a Cauchy argument.
\begin{theorem}\label{cauchy}
Suppose that $\var^0\in L^\infty (0,T)$.
The sequence $(\varepsilon^k)_{k\in\bbN}$ defined by
$(\ref{algo1})-(\ref{algo4})$ is convergent in $L^2(0,T)$.
\end{theorem}
{\bf Proof}: We still denote by $l_{\var^0}$ the value of $J$ on $C_{\var^0}$
and by $\widetilde{J}$ the shifted cost functional $J-l_{\var^0}$.
Suppose first that $\forall k \in
\bbN,\ \widetilde{J}(\varepsilon^{k})\neq 0$. By (\ref{dist}),
there exists $k_0$ such that (\ref{estloja}) holds (with
$\widetilde{C}_J=C_{\varepsilon^0}$) for all $\varepsilon^k$ with
$k\geq k_0$. Consider an integer $k\geq k_0$. We have:  
\beqn 
\left(\big(\widetilde{J}(\varepsilon^{k})\big)^{\wt}-\big(\widetilde{J}(\varepsilon^{k+1})\big)^{\wt}\right)
&\geq & \frac{\wt}{(\widetilde{J}(\varepsilon^{k+1}))^{1-\wt}}
\big(J(\varepsilon^{k+1})-J(\varepsilon^k)\big)\label{in1}\\
&\geq&\frac{\wk\wt}{\|\nabla J(\varepsilon^{k+1})\|_{L^1(0,T)}}
\Big(\Big(\frac{2}{\delta}-1\Big) \| \varepsilon^{k+1}- \tilde{\var}^k \|_{L^2(0,T)}^2\nonumber\\
&&\phantom{\frac{\wk\wt}{\|\nabla J(\varepsilon^{k+1})\|_{L^1(0,T)}}
\Big(} +\Big(\frac{2}{\eta}-1 \Big) \| \widetilde{\varepsilon}^k - \varepsilon^k \|^2_{L^2(0,T)} \Big)
\label{in2}\\
&\geq&\frac{\wk\wt a_{(\delta,\eta)}}{\lambda}\big( \| \varepsilon^{k+1}-
\tilde{\var}^k \|_{L^2(0,T)}+\!\| \tilde{\var}^{k}- \varepsilon^k
\|_{L^2(0,T)}\big)\label{nlbnl}\\
&\geq&\frac{\wk\wt a_{(\delta,\eta)}}{\lambda}\|
\varepsilon^{k+1}- \varepsilon^k \|_{L^2(0,T)}, \nonumber \eeqn
where $\ds
a_{(\delta,\eta)}=\frac{1}{\max(\delta,\eta)}-\frac{1}{2}$. 
The inequality (\ref{in1}) comes from the concavity of $s\mapsto s^{\wt}$, whereas
(\ref{in2}) is a consequence of (\ref{estloja}) and (\ref{monotonic}). Inequality
(\ref{nlbnl}) follows from (\ref{estnab}).\\
Since
$\left(\left(\widetilde{J}(\varepsilon^{k})\right)^{\wt}\right)_{k\in\bbN}$
is a Cauchy sequence (as a monotonic sequence bounded by
$(2\|O\|_*)^{\wt}$), we obtain that $(\varepsilon^k)_{k\in\bbN}$ is also
a Cauchy sequence.\\ 
If there exists $k_1$ such that
$\widetilde{J}(\varepsilon^{k_1})=0$, the monotonicity of the
algorithm implies that
$$J(\varepsilon^{k_1})=J(\varepsilon^{k_1+1})=J(\varepsilon^{k_1+2})=...$$
and by (\ref{monotonic}) the sequence $(\varepsilon^k)_{k\in\bbN}$
is constant for $k\geq k_1$.

\qed

\begin{remark}{\bf :} Thanks to the definition of the sequence $(\varepsilon^k)_{k\in\bbN}$ and to the regularity of the solutions $\psi$ and $\chi$ of the appropriate Schr\"odinger equations (see lemmas $\ref{schrod0}$ to $\ref{schrod2}$), we can easily prove by induction that if $\var^0 \in W^{1,\infty}(0,T)$, then for all $k\in\bbN$, 
$\var^k \in W^{1,\infty}(0,T)$.
\end{remark}

\section{Rate of convergence}\label{sec8}
The rate of convergence can be now evaluated by a second use of the
\L ojasiewicz inequality. The result is summarized in the next theorem.
\begin{theorem}
Let us denote by $\var^\infty$, the limit of $(\var^k)_{k\in\bbN}$ defined by
$(\ref{algo1})-(\ref{algo2})$ and $\wt$, $\wk$ the real numbers appearing
in $(\ref{estloja})$, where $C_{\var^0}=\{\var^\infty\}$.\\
\indent If $\wt<\dfrac{1}{2}$, then there exists $c>0$ such that
$\|\var^k-\var^\infty\|_{L^2(0,T)}\leq c k^{-\frac{\tit}{1-2\tit}}$.\\
\indent If $\wt=\dfrac{1}{2}$, then there exist $c'$ and $\tau$ such that:
\beq\label{optim}
\|\var^k-\var^\infty\|_{L^2(0,T)}\leq c'e^{-\tau k}.
\eeq

\end{theorem}
{\bf Proof:}\indent As in the proof of Theorem \ref{cauchy}, let be $k_0$, an integer such
that 
\beq\label{loj}
\forall \ell\geq k_0  \qquad \|\nabla J (\var^\ell)\|_{L^1(0,T)}\geq\wk |\widetilde{J}(\var^\ell)|^{1-\wt}.\eeq
Let us fix $k\geq k_0$ and introduce $\Delta^k$ defined by:
\beq\nonumber
\Delta^k=\sum_{\ell=k}^{\infty}\|\var^{\ell+1}-\tilde{\var}^\ell\|_{L^2(0,T)}+\|\tilde{\var}^{\ell}-\var^\ell\|_{L^2(0,T)}.
\eeq
With no loss of generality, we may assume that $\Delta^k>0$ for all
$k\geq k_0$.
Summing (\ref{nlbnl}) between $k$ and $+\infty$, we obtain:
\beq\nonumber
\big(\widetilde{J}(\var^{k})\big)^{\wt}\geq \frac{\wk\tit a_{(\delta,\eta)}}{\lambda}\Delta^k.
\eeq
This estimate, combined with (\ref{loj}), with $\ell=k$ yields:
\beq\nonumber
 \|\nabla J (\var^k)\|_{L^1(0,T)}
\geq \wk\Big(\frac{\wk\wt
  a_{(\delta,\eta)}}{\lambda}\Delta^k\Big)^{\frac{1-\tit}{\tit}}.
\eeq
From Lemma~\ref{estim}, we obtain:
\beq\nonumber
 \lambda (\Delta^{k-1}-\Delta^k)
\geq \wk\Big(\frac{\wk\wt
  a_{(\delta,\eta)}}{\lambda}\Delta^k\Big)^{\frac{1-\tit}{\tit}},
\eeq
which may be written as follows:
\beq\label{gene}
\frac{\Delta^{k-1}-\Delta^k}{(\Delta^k)^\beta}
\geq \upsilon,
\eeq
with $\beta=\dfrac{1-\tit}{\tit}$ and $\upsilon=\frac{\wk}{\lambda}\Big(\frac{\wk\wt
  a_{(\delta,\eta)}}{\lambda}\Big)^{\frac{1-\tit}{\tit}}$. Suppose now that
 $\wt=\dfrac{1}{2}$,{\it i.e.}, $\beta=1$. The equation (\ref{gene}) then
 becomes:
\beq\nonumber
(1+\upsilon)^{k_0}\Delta^{k_0}\big(\frac{1}{1+\upsilon}\big)^{k}\geq \Delta^k,
\eeq
and (\ref{optim}) is proved with $c'=(1+\upsilon)^{k_0}\Delta^{k_0}$ and $\tau=\ln(1+\upsilon)$.\\ Suppose now that $\wt<\frac{1}{2}$.
Let be $r\in]0,1[$, and suppose first that:
\beq\nonumber
(\Delta^k)^\beta \geq r(\Delta^{k-1})^\beta.
\eeq
Since $1-\beta<0$, the function $s\mapsto s^{1-\beta}$ is concave
 and we have:
\beq
(\Delta^{k})^{1-\beta}-(\Delta^{k-1})^{1-\beta}
\geq(\beta-1)\frac{\Delta^{k-1}-\Delta^k}{(\Delta^{k-1})^\beta}
\geq (\beta-1)r\frac{\Delta^{k-1}-\Delta^k}{(\Delta^{k})^\beta}
\geq (\beta-1)r\upsilon.\nonumber
\eeq
In the other case:
\beq
(\Delta^{k})^{1-\beta}-(\Delta^{k-1})^{1-\beta}\geq(\Delta^{k})^{1-\beta}-(r^{\frac{1}{\beta}}\Delta^{k})^{1-\beta}
=(1-r^{\frac{1-\beta}{\beta}})(\Delta^k)^{1-\beta}
\geq (1-r^{\frac{1-\beta}{\beta}})(\Delta^{k_0})^{1-\beta}.\nonumber
\eeq
Thus, in any case, there exists $\upsilon'>0$ independent of $k$, such that:
\beq\label{ineq}
(\Delta^{k})^{1-\beta}-(\Delta^{k-1})^{1-\beta}\geq\upsilon'.
\eeq
Consider now $k'>k$, the inequality (\ref{ineq}) implies that for a
small enough $c$, one have:
\beq\nonumber
\Delta^{k'}\leq
\Big(\upsilon'(k'-k)+(\Delta^{k})^{2-\frac{1}{\tit}}\Big)^{-\frac{\tit}{1-2\tit}}\leq
c k'^{-\frac{\tit}{1-2\tit}},
\eeq
and the result follows.

\qed

\begin{remark}{\bf:} Thanks to Lemma~$\ref{esthess}$, we have thus
  obtained that if $\alpha> 6T\|\mu\|^2_{L^\infty}\|O\|_*$ the convergence of
  the sequence is at least linear.
\end{remark}

{\bf Acknowledgment:} This work has been initiated at the summer school and
workshop ``Partial Differential Equations, Optimal Design and Numerics",
organized by  G. Buttazzo and E. Zuazua and is partially supported by the A.C.I
  "Simulation Mol\'eculaire'' of the french MENRT and by the Deutsche Forschungsgemeinschaft, SFB 404, B8. J.S. acknowledges helpful  discussions on that topic with J. Bolte (Laboratoire de Combinatoire et Optimisation, Universit\'e Pierre \& Marie  Curie, Paris) and O. Kavian (Laboratoire de Math\'ematiques , Universit\'e de Versailles Saint-Quentin). L.B. thanks G. Turinici  (CEREMADE, Universit\'e Paris-Dauphine).

\section*{Appendix: Proof of Lemma \ref{loja2}}
Consider $\var\in L^2(0,T)$ and $J$ defined by (\ref{defJ}). For reason of simplicity, we
suppose that $J(\var)=0$, $\nabla J(\var)=0$.\\
Thanks to Lemma \ref{hes}, the operator $H_J(\var)$ is a Fredholm
operator. The Fredholm alternative states then that either $H_J(\var)$ is
bijective or $Ker H_J(\var)=span(\varphi_1,...,\varphi_m)$, with $m>0$. Let us
denote by $\Pi$, the orthogonal
projection on $Ker H_J(\var)$ (with $\Pi=0$ if $Ker H_J(\var)=0$). The operator $L=\Pi+H_J(\var)$ is then bijective on
$L^2(0,T)$.\\
We are now in the position to apply the local inverse mapping theorem to
$\mathcal{L}=\Pi+\nabla J$ (analytic version, see \cite{zeidler}, Corollary 4.37, p.172),
that asserts there exist $V$ and $V'$ two neighborhoods of 0 in
$L^2(0,T)$ and $K:V'\rightarrow V$ an analytic mapping such that:
\beq
\forall \var'\in V,\ K(\mathcal{L}(\var'))=\var',\qquad 
\forall \var''\in V',\ \mathcal{L}(K(\var''))=\var''.\nonumber
\eeq
Since $\mathcal{L}$ and $K$ are $C^\infty$, there exist $C$ and $C'$ such that:
\beqn
\forall \var_1,\var_2\in V,\qquad
\|\mathcal{L}(\var_2)-\mathcal{L}(\var_1)\|_{L^2(0,T)}&\leq &C \|\var_2-\var_1\|_{L^2(0,T)}\nonumber\\
\forall \var'_1,\var'_2\in V',\qquad\|K(\var'_2)-K(\var'_1)\|_{L^2(0,T)}&\leq &C' \|\var'_2-\var'_1\|_{L^2(0,T)}.\nonumber
\eeqn 
Consider now $\var'\in V\cap V'$.
For $\zeta \in \bbR^m$ such that $\sum_{j=1}^m\zeta_j\varphi_j\in V$, let us define
$\Gamma:\zeta\mapsto J\big(K(\sum_{j=1}^m\zeta_j\varphi_j)\big),$
and $\xi \in \bbR^m$ such that $\Pi \var'=\sum_{j=1}^m \xi_j \varphi_j$. Let us
first estimate $\nabla\Gamma (\xi)$. Using $\Pi \var'\in V'$, we obtain:
\beqn
|\nabla\Gamma (\xi)|&\leq&C'' \|\nabla J\big(K(\Pi \var')\big)\|_{L^2(0,T)}
=C'' \|\nabla
J(\var')+\nabla J\big(K(\Pi \var')\big)-\nabla J(\var')\|_{L^2(0,T)}\nonumber\\
&\leq&C''\big(\|\nabla
J(\var')\|_{L^2(0,T)}+C \|K(\Pi \var')-\var'\|_{L^2(0,T)}\big)\nonumber\\
&=&C''\big(\|\nabla
J(\var')\|_{L^2(0,T)}+C \|K(\Pi \var')-K\big(\Pi\var'+\nabla J
(\var')\big)\|_{L^2(0,T)}\big)\nonumber\\
&\leq& c\|\nabla J(\var')\|_{L^2(0,T)},\label{eq1}
\eeqn 
where $c=C''(1+C C')$. On the other hand, one has:
\beqn
|J(\var')-\Gamma (\xi)|&=&|J(\var')-J\big(K(\Pi \var')\big)|
=\left|\int_0^1\frac{d}{d s} J\Big(\var'+s \big(K(\Pi
\var')-\var'\big)\Big)ds\right|\nonumber\\
&=&\left|\int_0^1 \left(\nabla J\Big(\var'+s \big(K(\Pi \var')-\var'\big)\Big),K(\Pi
\var')-\var'\right)ds\right|\nonumber\\
&\leq&\|K(\Pi \var')-\var')\|_{L^2(0,T)}\int_0^1\|\nabla J(\var')\|_{L^2(0,T)}+C s \|K(\Pi
\var')-\var'\|_{L^2(0,T)}ds \nonumber\\
&=&\|K(\Pi \var')-\var')\|_{L^2(0,T)} \Big(\|\nabla J(\var')\|_{L^2(0,T)}+\frac{C}{2} \|K(\Pi
\var')-\var'\|_{L^2(0,T)}\Big) \nonumber\\
 &\leq& c'\|\nabla J (\var')\|^2_{L^2(0,T)},\label{eq2}
\eeqn
where $c'=C'(1+\frac{CC'}{2})$.
By diminishing $V$, the \L ojasiewicz inequality (\ref{eqloja}) applied to the analytic functional $\Gamma$
states that there exist $ \theta\in ]0,1/2], \ \sigma>0$ such that:
\beqn
|\nabla \Gamma (\xi)|&\geq&|\Gamma (\xi)|^{1-\theta}
=|J(\var')-\Gamma (\xi)-J(\var')|^{1-\theta}\nonumber\\
&\geq&\frac{1}{2}|J(\var')|^{1-\theta}-\frac{1}{2}|J(\var')-\Gamma (\xi)|^{1-\theta}.\nonumber
\eeqn
Combining (\ref{eq1}) and (\ref{eq2}), we obtain:
\beqn
c\|\nabla J(\var')\|_{L^2(0,T)}&\geq&\frac{1}{2}|J(\var')|^{1-\theta}-c'\|\nabla J (\var')\|_{L^2(0,T)}^{2(1-\theta)},\nonumber
\eeqn
and the result follows.


\small
\begin{thebibliography}{99}
\bibitem{assion}A.~Assion, T.~Baumert, M.~Bergt, T.~Brixner,
  B.~Kiefer, V.~Seyfried, M.~Strehle, G.~Gerber, \emph{Control of
    chemical reactions by feedback-optimized phase-shaped femtosecond
    laser pulses}, Science, 282 (1998) 919--922.
\bibitem{B-M-S} J. M.~Ball, J. E.~Marsden, M.~Slemrod,
\emph{Controlability for distributed bilinear systems},
SIAM J. Cont. Opt., 20 (1982) 575--597.
\bibitem{baudouin} L.~Baudouin, O.~Kavian, J.-P.~Puel,
\emph{Regularity for a Schr\"odinger equation with singular potentials and
application to bilinear optimal control},
J. Diff. Eq., 216 (2005) 188--222.
\bibitem{baudouin2}
L.~Baudouin, 
\emph{A bilinear optimal control problem applied to a time dependent
Hartree-Fock equation coupled with classical nuclear dynamics},  
Portugaliae Mathematica (N.S.), 63 (1)  (2006) 293--325.  
\bibitem{baudouin3}
L.~Baudouin, 
\emph{Existence and regularity of the solution of a time dependent
Hartree-Fock equation coupled with a classical nuclear dynamics}, 
Revista Matematica Complutense, 18 (2) (2005) 285--314.  
\bibitem{beauchard}K.~Beauchard,
\emph{ Local controllability of a 1D Schr\"odinger equation }, J. Math. Pures Appl., 84 (7) (2005) 851--956.
\bibitem{brixner1} T.~Brixner, N. H.~Damrauer, P.~Niklaus, G.~Gerber, 
\emph{Photoselective adaptive femtosecond quantum control in the liquid phase}, 
Nature, 414 (2001) 57--60.
\bibitem{brown} E.~Brown, H.~Rabitz, 
\emph{Some mathematical and algorithmic challenges in the control of quantum dynamics phenomena}, 
J. Math. Chem., 31 (2002) 17--63.
\bibitem{pilot} E.~Canc\`es, C.~Le Bris, M.~Pilot,  
\emph{Optimal bilinear control for a Schr\"odinger equation}, 
C. R. Acad. Sci. Paris,  330 (S\'erie 1)  (2000) 567--571.
\bibitem{bookgt} E.~Canc\`es, C.~Le Bris, Y.~Maday, G.~Turinici, 
\emph{Mathematical Foundations of Molecular Modelling}, 
Oxford Univ. Press, oxford, 2007.
\bibitem{Caze}  T.~Cazenave,
\emph{An introduction to nonlinear Schr\"odinger equation},
 third edition, Textos de M\'etodos Matem\'aticos 26, Rio de Janeiro, 1996.
\bibitem{kasparian} J.~Kasparian, M.~Rodriguez, G.~M\'ejean, J.~Yu,
  E.~Salmon, H.~Wille, R.~Bourayou, S.~Frey, Y.-B.~Andr,
  A.~Mysyrowicz, R.~Sauerbrey, J.-P.~Wolf, L.~Woste, 
\emph{White-light filaments for atmospheric analysis}, Science, 301 (2003) 61--64.
\bibitem{kavian}
A.~Haraux, M.A.~Jendoubi, O.~Kavian, 
\emph{Rate of decay to equilibrium in some semilinear parabolic equations},
J. Evol. Equ., 3 (2003) 463--484.
\bibitem{I-K}
K.~Ito, K.~Kunisch, 
\emph{Optimal bilinear control of an abstract Schr{\"o}dinger equation}, 
SIAM J. Cont. Opt. , to appear (2007).
\bibitem{jendoubi}
M.~A.~Jendoubi, 
\emph{A simple unified approach to some convergence theorems of L. Simon},  
J. Func. Anal., 153 (1998) 187--202.
\bibitem{lebris}
C. Le Bris, \emph{Computational Chemistry, Handbook of Numerical Analysis}, Ph. G. Ciarlet ed., volume X, North-Holland, 2003.  
\bibitem{loja}
S.~{\L}ojasiewicz, 
\emph{Une propri\'et\'e topologique des sous-ensembles analytiques r\'eels}, 
Colloques internationaux du CNRS. Les \'equations aux d\'eriv\'ees partielles, 117 (1963).
\bibitem{loja2}
S.~{\L}ojasiewicz, 
\emph{Sur la g\'eom\'etrie semi- et sous-analytique}, 
Ann. Inst. Fourier, 43 (1993) 1575--1595.
\bibitem{salomon}
Y.~Maday, J.~Salomon, G.~Turinici, 
\emph{Monotonic time-discretized schemes in  quantum control}, 
Num. Math., 103 (2) (2006) 323--338.
\bibitem{maday}
Y.~Maday, G.~Turinici, 
\emph{New formulations of monotonically convergent quantum control algorithms},
J. Chem. Phys, 118 (18) (2003).
\bibitem{rabitz2} 
H.~Rabitz, R. de~Vivie-Riedle, M.~Motzkus, K.~Kompa, 
\emph{Whither the future of control-
         ling quantum phenomena?}, Science, 288 (2000) 824--828.
\bibitem{rabitz3} 
H.~Rabitz, G.~Turinici, E.~Brown, 
\emph{Control of quantum dynamics: Concepts, procedures and future prospects}, 
In Ph. G. Ciarlet, editor, Computational Chemistry, Special Volume
(C. Le Bris Editor) of Handbook of Numerical Analysis, vol X, Elsevier Science B.V., 2003.
\bibitem{R-S}
M. Reed, B. Simon,
\emph{Methods of Modern Mathematical Physics, II, Fourier analysis, self-adjointness},
Academic Press, 1975.
\bibitem{salomon-cv}
J.~Salomon, 
\emph{Convergence of the time-discretized monotonic schemes}, 
M2AN, 41 (1) (2007) 77--93.
\bibitem{salomon-cdc05}
J.~Salomon, 
\emph{Limit points of the monotonic schemes in quantum control}, 
Proceedings of the 44th IEEE Conference on Decision and Control, Sevilla, 2005.
\bibitem{these}
J.~Salomon, 
\emph{Contr\^ole en chimie quantique : conception et analyse de sch\'emas d'optimisation}, 
Th\`ese de l'Universit\'e Pierre et Marie Curie, 2005.
\bibitem{shi} 
S.~Shi, A.~Woody, H.~Rabitz, 
\emph{Optimal control of selective vibrational excitation in harmonic linear chain molecules}, 
J. Chem. Phys., 88 (1988) 6870--6883.
\bibitem{simon}
L.~Simon, 
\emph{Asymptotics for a class of non-linear evolution equations},
with applications to geometric problems, Ann. of Math., 118 (1983) 525--571.
\bibitem{tannor}
D.~Tannor, V.~Kazakov, V.~Orlov, 
\emph{Control of Photochemical Branching: Novel Procedures for Finding Optimal Pulses and Global Upper Bounds}, 
in {J.~Broeckhove}, {L.~Lathouwers } (Eds), Time Dependent Quantum Molecular Dynamics, Plenum Press, New York,  1992, 347--360.
\bibitem{tur_cdc} G.~Turinici, 
\emph{Monotonically Convergent Algorithms for Bounded Quantum Controls}, 
Proceedings of the LHMNLC03 IFAC Conference, 2003 263--266.
\bibitem{vogt} 
G. Vogt, G. Krampert, P. Niklaus, P. Nuernberger, G. Gerber, 
\emph{Optimal control of photoisomerization}, 
Phys. Rev. Lett., 94: 68305 (2005).
\bibitem{weinacht} T.~Weinacht, J.~Ahn, P.~Bucksbaum, \emph{Controlling the shape of a quantum wavefunction},
Nature, 397 (1999) 233--235.
\bibitem{zeidler}
E.~Zeidler, 
\emph{Nonlinear functional analysis and its applications},
tome 1, Springer-Verlag, Berlin/New York, 1985.
\bibitem{zhu}
W.~Zhu, H.~Rabitz, 
\emph{A rapid monotonically convergent algorithm for quantum optimal control over the expectation value of a definite operator}, 
J. Chem. Phys., 109 (1998) 385--391.
\end{thebibliography}
\end{document}